\documentclass[12pt]{amsart}

\usepackage{amsmath, amscd, amssymb, latexsym}
\usepackage{graphicx}
\usepackage{hyperref}
\usepackage{cleveref}
\usepackage{mathptmx}
\usepackage{esint}
\usepackage{booktabs}
\usepackage{tikz} 

\theoremstyle{plain}

\newtheorem{theorem}{Theorem}[section]
\newtheorem{lemma}[theorem]{Lemma}

\newtheorem{corollary}[theorem]{Corollary}

\newtheorem{question}[theorem]{Question}

\newtheorem{theoremx}{Theorem}

\newcommand{\CP}{\mathbb{CP}}
\newcommand{\W}{\mathcal W}
\newcommand{\Sc}{\mathcal S}
\newcommand{\Yadm}{\mathcal Y_2^+}

\begin{document}

\title[Normalized Weyl Energy]{Normalized Weyl Energy and the\\
Differential Topology of Four-Manifolds}
\author{Mijia Lai}
\address{Shanghai Jiao Tong University, Shanghai, 200240, China}
\email{laimijia@sjtu.edu}
\date{}

\begin{abstract}
The Weyl functional $\int_M\|W_g\|^2\,d\mu_g$ plays an important role in four dimensional conformal geometry.
Using this functional, Chang,
Gursky, and Zhang introduced an invariant of an admissible conformal
class and a smooth invariant of the underlying manifold obtained by
taking the infimum over all such classes. We use a normalization in
which both invariants, denoted by $\beta(M,c)$ and $\beta(M)$,
take values in $[0,1)$.

We prove that every admissible four-manifold $M$ with $\beta(M)$
sufficiently close to $1/2$ is diffeomorphic to $\mathbb{CP}^2$ up
to reversal of orientation. At $\beta(M)=1/2$, this result and the
conformal characterization of the Fubini--Study class confirm
Conjecture~1.1 of Chang, Gursky, and Zhang. Among admissible
four-manifolds with indefinite intersection form, we prove the
corresponding stability of the diffeomorphism type of $S^2\times S^2$
near $2/3$. For the one- and two-point blow-ups of $\mathbb{CP}^2$,
we show that $\beta(M)$ lies strictly above the respective
topological lower bounds of $2/3$ and $23/30$.

These results lead us to ask whether $1$ is the only accumulation
point of the values of $\beta(M)$ and whether, among simply connected
admissible four-manifolds, $\beta(M)$ determines the diffeomorphism
type up to reversal of orientation.
\end{abstract}

\maketitle


\section{Introduction}
\label{sec:introduction}

On a closed, connected, orientable surface $\Sigma$, a single curvature
integral determines the diffeomorphism type. Indeed, the Gauss--Bonnet
formula gives
\[
   \int_\Sigma K_g\,dA_g
      =2\pi\chi(\Sigma).
\]
Uniformization gives this identity a geometric interpretation: according
to the sign of $\chi(\Sigma)$, every conformal class contains a metric
of constant positive, zero, or negative curvature. In dimension three,
Thurston's geometrization program, completed by Perelman using
Hamilton's Ricci flow, describes the topology of a closed orientable
three-manifold through its decomposition into pieces modeled on the
eight Thurston geometries \cite{MorganTian}. 

In dimensions four and higher, it is considerably more difficult to
deduce a topological classification from curvature conditions.
In dimension four, however, the intersection form provides a link
between curvature and topology. For a closed, simply connected
smooth four-manifold, Freedman's theorem shows that the integral
intersection form determines the oriented homeomorphism type, while
Donaldson's theorem shows that a definite intersection form must be
diagonalizable over the integers \cite{Freedman,Donaldson}. From the
geometric side, the Weyl energy $\int_M\|W_g\|^2\,d\mu_g$ is
conformally invariant in dimension four. Its self-dual and
anti-self-dual parts enter the Hirzebruch signature formula  and the Chern--Gauss--Bonnet formula. The former identity gives the signature of the intersection
form, while the latter gives its rank through
$b_2(M)=\chi(M)-2$ in the simply connected case. These identities
allow integral curvature inequalities to constrain the intersection
form and, under suitable pinching hypotheses, to determine the
diffeomorphism type of the underlying manifold, see e.g.,
\cite{CGYsphere,CGZ}.

Let $M$ be a closed, connected, oriented smooth four-manifold.
For a Riemannian metric $g$ on $M$, write
$E_g=\operatorname{Ric}_g-\tfrac14 R_g g$ for its trace-free
Ricci tensor and $W_g^\pm$ for its Weyl curvature operators on
$\Lambda^\pm$. Using the Hilbert--Schmidt norm on these
operators, set
\[
    \mathcal W^\pm([g])
       =\int_M\|W_g^\pm\|^2\,d\mu_g,
    \qquad
    \mathcal W([g])=\mathcal W^+([g])+\mathcal W^-([g]),
\]
and define
\[
    \mathcal S([g])
       =\int_M\left(\frac{R_g^2}{24}
                    -\frac{|E_g|^2}{2}\right)d\mu_g.
\]
The Chern--Gauss--Bonnet and signature formulas then read
\begin{equation}\label{eq:characteristic}
   8\pi^2\chi(M)=\W([g])+\Sc([g]),\qquad
   12\pi^2\tau(M)=\W^+([g])-\W^-([g]).
\end{equation}
Both Weyl energies are invariant under conformal changes of metric in
dimension four, so the first identity gives the same invariance for
$\Sc$. In terms of the Schouten tensor
$P_g=\tfrac12(\operatorname{Ric}_g-\tfrac16R_g g)$, one has
$\Sc([g])=4\int_M\sigma_2(g^{-1}P_g)\,d\mu_g$.

Following Chang, Gursky, and Zhang \cite{CGZ}, let $\Yadm(M)$ denote
the set of conformal classes $c$ with positive Yamabe constant $Y(M,c)$
and $\Sc(c)>0$. We call such a class \emph{admissible}, and call $M$
admissible when $\Yadm(M)\neq\varnothing$. An admissible manifold
necessarily has positive Euler characteristic in view of  \eqref{eq:characteristic}. We define the
\emph{normalized Weyl energy} of a conformal class as 
\begin{equation}\label{eq:beta}
   \beta(M,c)=\frac{\W(c)}{8\pi^2\chi(M)}
   \end{equation}
 Taking the infimum 
 \begin{equation}
   \beta(M)=\inf_{c\in\Yadm(M)}\beta(M,c).
\end{equation}
gives a smooth invariant, which we call \emph{the normalized Weyl energy} of $M$. These two invariants are
first introduced by Chang, Gursky and Zhang \cite{CGZ}. We slightly modify the normalization so that $\beta(M)\in[0,1)$. If $B$ denotes the invariant in the original
normalization of Chang-Gursky-Zhang\cite{CGZ}, then
$\beta=B/(B+4)$ for both a conformal class and a manifold. In
particular, their values $4$ and $8$ correspond to $1/2$ and $2/3$.

The remarkable conformal four sphere theorem of Chang,
Gursky, and Yang \cite{CGYsphere} implies that if an admissible
manifold satisfies $\beta(M)<1/2$, then $M$ is diffeomorphic to the
standard $S^4$ and consequently $\beta(M)=0$. 

Gursky's inequality \cite{Gursky} provides a lower bound of the normalized Weyl energy on admissible four-manifolds in terms of the signature and Euler characteristic number. 
If 
$b_2^+(M)>0$ and $Y(M,c)>0$, then
\begin{equation}\label{eq:gursky}
   \W^+(c)\geq\frac{4\pi^2}{3}
                   \bigl(2\chi(M)+3\tau(M)\bigr),
\end{equation}
with equality precisely when $c$ contains a positive
K\"ahler--Einstein metric compatible with the chosen orientation.
Since $\W=2\W^+-12\pi^2\tau(M)$, it follows that
\begin{equation}\label{eq:lower}
   \beta(M)\geq\frac23-\frac{\tau(M)}{2\chi(M)}
   \qquad\text{if } b_2^+(M)>0.
\end{equation}
When the intersection form is indefinite, applying the same estimate
in both orientations gives
\begin{equation}\label{eq:indef}
   \beta(M)\geq\frac23+\frac{|\tau(M)|}{2\chi(M)} .
\end{equation}

These inequalities give exact values whenever an appropriate equality
metric is available. For example, the Fubini--Study metric $\omega_{FS}$ on $\CP^2$ and the
standard product metric $g_{rd}\oplus g_{rd}$ on $S^2\times S^2$ realize the values
$1/2$ and $2/3$, respectively. More generally, write
$X_k=\CP^2\#k\overline{\CP}^{\,2}$ with its standard orientation, then $\chi(X_k)=3+k$ and $\tau(X_k)=1-k$. The positive
K\"ahler--Einstein metrics on del Pezzo surfaces given by Tian-Yau and Tian
\cite{TianYau,Tian} realize Gursky's lower bound and give
\begin{equation}\label{eq:delpezzo}
   \beta(X_k)=\frac{7k+9}{6(k+3)}
   \qquad\text{for } k=3,\ldots,8 .
\end{equation}
A different source of exact values is the signature inequality
$\W\geq12\pi^2|\tau|$. The self-dual metrics of positive Ricci
curvature constructed by LeBrun, Nayatani, and Nitta \cite{LNN}
realize this bound on $\#_n\CP^2$ for $n=2,3$. Their classes are
admissible since $\Sc=4\pi^2(4-n)>0$.

\begin{table}[htbp]
\centering
\begin{tabular}{@{}lc@{\qquad\qquad}lc@{}}
\toprule
$M$ & $\beta(M)$ & $M$ & $\beta(M)$\\
\midrule
$S^4$                                  & $0$    & $X_3$ & $5/6$\\
$\CP^2$                                & $1/2$  & $X_4$ & $37/42$\\
$S^2\times S^2$                        & $2/3$  & $X_5$ & $11/12$\\
$(S^2\times S^2)/\langle\iota\rangle$  & $2/3$  & $X_6$ & $17/18$\\
$\CP^2\#\CP^2$                         & $3/4$  & $X_7$ & $29/30$\\
$\CP^2\#\CP^2\#\CP^2$                  & $9/10$ & $X_8$ & $65/66$\\
\bottomrule
\end{tabular}
\caption{Selected exact values of the normalized Weyl energy.
Here $\iota(x,y)=(-x,-y)$ is the diagonal antipodal involution.}
\label{tab:values}
\end{table}

The results of Chang, Gursky, and Zhang \cite{CGZ} exhibit both stability of the
diffeomorphism type under pinching of the normalized Weyl energy and
conformal rigidity at equality. They
proved that if $b_2^+(M)>0$ and some admissible conformal class $c$
satisfies $1/2\leq\beta(M,c)<1/2+\delta$ for a universal $\delta>0$,
then $M$ is diffeomorphic to $\CP^2$ \cite[Theorem~A]{CGZ}.
Under the same Betti number hypothesis, if the manifold invariant
satisfies $\beta(M)=1/2$, then the minimizing conformal classes are
precisely the pullbacks of the $[\omega_{FS}]$ by diffeomorphisms 
from $M$ to $\CP^2$ \cite[Theorem~B]{CGZ}. 

At the product value, Chang, Gursky, and Zhang proved that if $M$ has
indefinite intersection form and an admissible class $c$ satisfies
$\beta(M,c)=2/3$, then $(M,c)$ is conformally diffeomorphic to
$(S^2\times S^2, [g_{rd}\oplus g_{rd}]) $ \cite[Theorem~C]{CGZ}. The diagonal antipodal quotient
in Table~\ref{tab:values} explains the hypothesis on the intersection
form. It has $b_2=0$ and fundamental group $\mathbb{Z}/2$, yet its
invariant is also $2/3$: the round product metric descends to the
quotient, and every admissible class on the quotient lifts to the
product with the same normalized Weyl energy. Within the class of
Bach-flat metrics, Zhang extended this conformal rigidity to a
neighborhood of the product value. If $b_2^+(M)=b_2^-(M)>0$ and
$c$ is an admissible Bach-flat class satisfying
$2/3\leq\beta(M,c)<2/3+\delta'$ for a universal $\delta'>0$, the
same conformal classification holds \cite[Theorem~B]{Zhang}.
Zhang's proof uses Yamabe representatives with small $L^2$ norm of
the trace-free Ricci tensor and a lower bound for their Yamabe
constants that excludes nonflat Ricci-flat bubbles
\cite[\S4.2]{Zhang}. To establish stability of the diffeomorphism
type under pinching of $\beta(M)$, one must instead obtain compactness
for a minimizing sequence of admissible conformal classes, whose
members need not be Bach-flat.

Based on these results, we continue their study in terms of the 
normalized Weyl energy and their relationship to the
differential topology of $M$. 

\begin{question}
\begin{enumerate}
\item[(a)] Is $1/2$ an isolated value of $\beta$ among closed,
oriented, admissible four-manifolds? If $\beta(M)=1/2$, must $M$
be diffeomorphic to $\mathbb{CP}^2$ up to reversal of orientation? (Conjecture 1.1 of \cite{CGZ})

\item[(b)] Among closed, oriented, admissible four-manifolds with
indefinite intersection form, is $2/3$ an isolated value of $\beta$?
If $\beta(M)=2/3$ for such a manifold $M$, must $M$ be
diffeomorphic to $S^2\times S^2$?

\end{enumerate}
\end{question}

A broader question could be   
\begin{question}
 If $M$ and $N$ are simply connected admissible
four-manifolds with $\beta(M)=\beta(N)$, must they be
diffeomorphic, possibly after reversing the orientation of $N$?
\end{question}

By Table~\ref{tab:values}, the simply-connectedness is necessary. 
At the value $2/3$, it asks whether one can remove 
the indefinite assumption on the intersection form of $M$. 
It leads naturally to the possibility of exotic smooth structures on $S^4$ and
$\mathbb{CP}^2$. Let $M$ be a simply connected admissible
four-manifold. If its intersection form is indefinite, Theorem~B
identifies $M$ with $S^2\times S^2$ whenever $\beta(M)$ is
sufficiently close to $2/3$. Otherwise, $|\tau(M)|=b_2(M)$ and
$\chi(M)=2+b_2(M)$, so the signature inequality gives
\[
   \beta(M)\geq \frac{3|\tau(M)|}{2\chi(M)}
   =\frac{3b_2(M)}{2(b_2(M)+2)}.
\]
For $b_2(M)\geq2$, the right-hand side is at least $3/4$.
Thus only $b_2(M)=0$ and $b_2(M)=1$ remain near $2/3$.
By Freedman's classification \cite{Freedman}, the corresponding
underlying topological manifolds are $S^4$ and $\mathbb{CP}^2$,
respectively, up to orientation. Their standard smooth structures
have $\beta$-values $0$ and $1/2$. Consequently, a simply connected
counterexample to uniqueness at $2/3$ would have to be an admissible
exotic $S^4$ or $\mathbb{CP}^2$ with $\beta=2/3$. 

The one- and two-point blow-ups $X_1$ and $X_2$ are the only
diffeomorphism types of del Pezzo surfaces missing from
Table~\ref{tab:values}. Neither del Pezzo surface admits a
K\"ahler--Einstein metric \cite{Tian}. Both nevertheless carry a
positive Einstein metric conformal to an extremal K\"ahler metric,
namely the Page metric on $X_1$ \cite{Page} and the
Chen--LeBrun--Weber metric on $X_2$ \cite{CLW}. So we may ask:

\begin{question}
\label{q:delpezzo}
What are the exact values of $\beta(X_1)$ and $\beta(X_2)$?
Are the respective infima realized by the Page conformal class
on $X_1$ and the Chen--LeBrun--Weber conformal class on $X_2$?
\end{question}

The minimization question has an antecedent in Koca's doctoral
thesis \cite[\S1.5, p.~17]{Koca}.  LeBrun subsequently proved the
minimizing property among conformal classes of symplectic type
invariant under the respective torus actions
\cite[Theorem~C]{LeBrunWeyl}.

Koca also observed that no individual conformal class of positive
Yamabe type on either manifold attains Gursky's lower bound.
This leaves open the possibility that a sequence of admissible
classes approaches the bound. Theorem~\ref{thm:main-blowups} (in Section 3)
excludes this possibility on both manifolds. Whether the Page
and Chen--LeBrun--Weber classes realize the respective infima
remains open.

The examples in Table~\ref{tab:values} raise the question of how
the values of the normalized Weyl energy are distributed in
$[0,1)$. One accumulation point can already be identified.
For $k\geq1$, let $M_k=\#_k(S^2\times S^2)$. The standard
product metric for $k=1$ and the recent constructions of positive
Einstein metrics for $k\geq2$ show that each $M_k$ is admissible
\cite{Ozuch,LiSun,BroderLafuente,ChuStern}. For any admissible
conformal class $c$ on a closed four-manifold $M$, evaluation on
a unit-volume Yamabe representative gives
\[
    8\pi^2\chi(M)\bigl(1-\beta(M,c)\bigr)
      =\mathcal S(c)
      \leq \frac{Y(M,c)^2}{24}
      \leq 16\pi^2.
\]
It follows that $\beta(M)\geq 1-2/\chi(M)$. Since
$\chi(M_k)=2(k+1)$, we obtain
\[
    \frac{k}{k+1}<\beta(M_k)<1. \qquad (\text{the former inequality is strict})
\]
Thus $\beta(M_k)\to1$, and $1$ is an accumulation point (see Figure\ref{axis} below). 

\begin{figure}[htbp]
  \centering
\begin{tikzpicture}[x=13cm,y=1cm, line cap=round, font=\normalsize]

  \coordinate (bound) at (0.812,0);
  \coordinate (value) at (0.867,0);
  \draw[line width=0.55pt] (0,0) -- (1,0);
  \foreach \pos/\model/\modelclass/\level in {
      0/{S^4}/{[g_{\mathrm{rd}}]}/{0},
      0.5/{\mathbb{CP}^2}/{[\omega_{\mathrm{FS}}]}/{\frac{1}{2}},
      0.666667/{S^2\!\times\! S^2}/{[g_{\mathrm{rd}}\oplus g_{\mathrm{rd}}]}/{\frac{2}{3}}} {
    \fill (\pos,0) circle (1.65pt);
    \node at (\pos,0.68) {$\model$};
    \node at (\pos,0.30) {$\scriptstyle\modelclass$};
    \node[below=4pt] at (\pos,0) {$\level$};
  }

  \draw[line width=0.55pt] (bound) ++(0,-0.065) -- ++(0,0.13);
  \node[below=4pt] at (bound) {$\frac{k}{k+1}$};
  \fill (value) circle (1.65pt);
  \node[above=4pt] at (value) {$M_k$};
  \node at (0.766,0.31) {$\cdots$};
  \node at (0.925,0.31) {$\cdots$};
  \draw[->,>=stealth,line width=0.5pt] (0.946,0.31) -- (0.995,0.31);

  \draw[fill=white,line width=0.65pt] (1,0) circle (2.1pt);
  \node[below=4pt] at (1,0) {$1$};
\end{tikzpicture}
    \caption{Selected values of the normalized Wely energy} \label{axis}
    \end{figure}
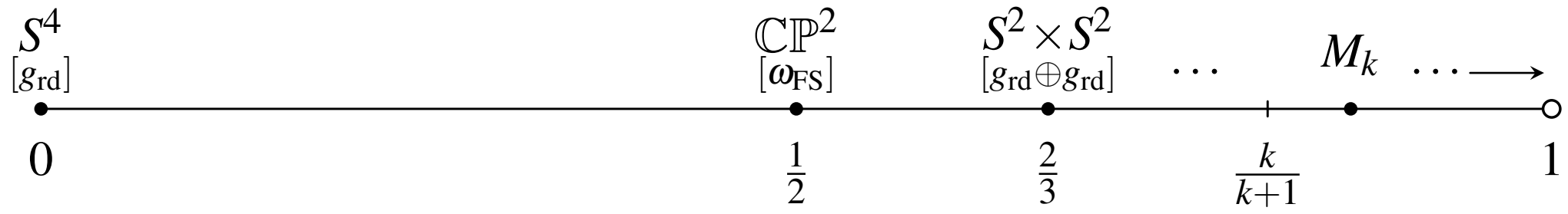
    
This leads to the following natural questions: 
\begin{question}\label{q:sublevel-finiteness}
Let $\mathcal B$ denote the set of values $\beta(M)$ as $M$
ranges over all closed, connected, oriented, admissible
four-manifolds.
\begin{enumerate}
    \item[a]  Is $\mathcal{B}\cap [0, \frac{2}{3})=\{0,\frac{1}{2}\}$? (Conjecture 1.2 of \cite{CGZ})
    \item[b]  Is $1$ the only accumulation point of $\mathcal B$?
\end{enumerate}

\end{question}

Our main results are as follows. In Section~\ref{sec:half-gap}, we
answer both parts of Question~1 affirmatively. Theorems~A and~B
establish stability of the diffeomorphism type under pinching of
$\beta(M)$. Theorem~A shows that every admissible four-manifold $M$
with $\beta(M)$ sufficiently close to $1/2$ is diffeomorphic to
$\mathbb{CP}^2$ up to reversal of orientation. At $\beta(M)=1/2$,
the equality characterization in \cite[Theorem~B]{CGZ} then shows
that the minimizing conformal classes are precisely the pullbacks of
the Fubini--Study class. This proves Conjecture~1.1 of \cite{CGZ}.
Theorem~B shows that, among admissible four-manifolds with indefinite
intersection form, every $M$ with $\beta(M)$ sufficiently close to
$2/3$ is diffeomorphic to $S^2\times S^2$.

In Section~\ref{sec:delpezzo-gaps}, we give a partial answer to
Question~\ref{q:delpezzo}. Theorem~C shows that the normalized Weyl
energies of the one- and two-point blow-ups of $\mathbb{CP}^2$ lie
strictly above their respective topological lower bounds:
\[
   \frac{2}{3}<\beta(X_1)\leq\beta(X_1,[g_{\mathrm{P}}]),
   \qquad
   \frac{23}{30}<\beta(X_2)\leq\beta(X_2,[g_{\mathrm{CLW}}]).
\]
Whether the Page and Chen--LeBrun--Weber conformal classes realize
the respective infima remains open.

\medskip
\noindent\textbf{Declaration of AI use.}
We used ChatGPT for proof exploration and editing.
We take responsibility for verification and exposition.

\section{Isolated values of the normalized Weyl energy}
\label{sec:half-gap}

In this section, we answer both parts of Question 1.1 affirmatively.

\begin{theoremx}
\label{thm:main-gap}
There exists a universal constant $\varepsilon_1\in (0, \frac{1}{2})$ such that,
for every closed, oriented, admissible four-manifold $M$,
the following are equivalent:
\begin{enumerate}
\item[(i)] $\bigl|\beta(M)-\tfrac12\bigr|<\varepsilon_1$.
\item[(ii)] $\beta(M)=\tfrac12$.
\item[(iii)] $M$ is diffeomorphic to $\mathbb{CP}^2$ up to
reversal of orientation.
\end{enumerate}
\end{theoremx}

\begin{theoremx}
\label{thm:main-product-gap}
There exists a universal constant $\varepsilon_2>0$ such that,
for every closed, oriented, admissible four-manifold $M$
with indefinite intersection form, the following are equivalent:
\begin{enumerate}
\item[(i)] $\bigl|\beta(M)-\tfrac23\bigr|<\varepsilon_2$.
\item[(ii)] $\beta(M)=\tfrac23$.
\item[(iii)] $M$ is diffeomorphic to $S^2\times S^2$.
\end{enumerate}
\end{theoremx}

Theorems~A and~B share a common compactness argument. We evolve the
unit-volume Yamabe representatives by the negative $L^2$-gradient
flow of the total squared Riemann curvature. Lemma~2.1 gives a
common existence interval and the scale-invariant curvature
estimate~(13). At each fixed positive time, Corollary~2.2 then
gives smooth convergence to a positive Einstein metric after
passing to a subsequence and pulling back by diffeomorphisms.

Throughout this section, write
\[
   Y_4=Y(S^4,[g_{\mathrm{rd}}])=8\sqrt6\,\pi,\qquad
   E_g=\operatorname{Ric}_g-\tfrac14R_g g,\qquad
   e(g)=\int_M |E_g|^2\,d\mu_g .
\]

To prove Theorem\ref{thm:main-gap}, it suffices to show $(i)\Longrightarrow (iii)$. Note $\beta(M)<\frac{1}{2}$ immediately forces $\beta=0$. We prove by contradiction, and thus may assume there are admissible manifolds $M_j$, none diffeomorphic to
$\CP^2$ in either orientation, and classes $c_j\in\Yadm(M_j)$
such that
\begin{equation}\label{eq:half-countersequence}
   \frac12\leq\beta(M_j)\leq\beta(M_j,c_j)
      <\frac12+\frac{2}{j}.
\end{equation}
Since $\beta(S^4)=0$, none of the $M_j$ is diffeomorphic to the
standard $S^4$. If $b_2(M_j)>0$, we can choose the orientation so
that $b_2^+(M_j)>0$. The gap theorem of Chang, Gursky, and Zhang
\cite[Theorem~A]{CGZ} then excludes all sufficiently large such
$j$. We may therefore assume $b_2(M_j)=0$.

An admissible class contains a metric of positive Ricci curvature
\cite{CGYricci}. Hence $\pi_1(M_j)$ is finite, $b_1(M_j)=0$,
$\chi(M_j)=2$, and $\tau(M_j)=0$. Let $d_j$ be the degree of the
universal cover $\widetilde M_j$. Since the cover is simply
connected and both characteristic numbers are multiplicative,
\[
   b_2^+(\widetilde M_j)+b_2^-(\widetilde M_j)=2d_j-2,
   \qquad
   b_2^+(\widetilde M_j)-b_2^-(\widetilde M_j)=0.
\]
Thus $b_2^+(\widetilde M_j)=b_2^-(\widetilde M_j)=d_j-1$.
If $d_j>1$, the cover has indefinite intersection form. The lifted
class is admissible and has the same normalized Weyl energy as
$c_j$, so \eqref{eq:indef} gives $\beta(M_j,c_j)\geq2/3$. This
contradicts \eqref{eq:half-countersequence} for large $j$.

Consequently $M_j$ is simply connected for large $j$. Integral
Poincar\'e duality and the universal coefficient theorem now give
$H_2(M_j;\mathbb Z)=0$. We write $\Sigma_j=M_j$ to indicate that
the remaining sequence consists of integral homology four-spheres.

Let $g_j\in c_j$ be a unit-volume Yamabe metric and put
$Y_j=Y(\Sigma_j,c_j)=R_{g_j}$.  At each point, let $\lambda_j$ be the
larger of the largest eigenvalues of $W^+_{g_j}$ and $W^-_{g_j}$.
The classification theorem of Chen and Zhu
\cite[Theorem~1.4]{ChenZhu} implies that a simply connected closed
four-manifold of positive Yamabe type satisfying
$36\int\lambda_j^2\,d\mu_{g_j}<Y_j^2$ is diffeomorphic to $S^4$.
So on $\Sigma_j$, we have conversely $Y_j^2\leq 36\int_{\Sigma_j}\lambda_j^2\,d\mu_{g_j}$. 
Moreover, a trace-free symmetric endomorphism $A$ of a
three-dimensional Euclidean space satisfies
$\lambda_{\max}(A)^2\leq\tfrac23\|A\|^2$. We conclude that
\begin{equation}\label{eq:half-chen-zhu}
   Y_j^2\leq36\int_{\Sigma_j}\lambda_j^2\,d\mu_{g_j}
          \leq24\W(c_j).
\end{equation}
Set $\eta_j=\beta(\Sigma_j,c_j)-1/2$. Since $\chi(\Sigma_j)=2$,
\eqref{eq:characteristic} gives
\[
   \W(c_j)=8\pi^2+16\pi^2\eta_j,\qquad
   \Sc(c_j)=\frac{Y_j^2}{24}-\frac12e(g_j)
           =8\pi^2-16\pi^2\eta_j .
\]
Combining these identities with \eqref{eq:half-chen-zhu}, we obtain
\begin{equation}\label{eq:half-initial}
   0\leq e(g_j)\leq64\pi^2\eta_j\longrightarrow0,\qquad
   Y_j\longrightarrow s_*:=8\sqrt3\,\pi=\frac{Y_4}{\sqrt2}.
\end{equation}

Next we list several basic facts about the negative $L^2$-gradient flow of the total curvature functional.
Let $(N,h)$ be a closed Riemannian four-manifold. Consider the
total curvature functional
\[
   \mathcal F(h)=\frac12\int_N|\operatorname{Rm}_h|^2\,d\mu_h,
\]
where $\operatorname{Rm}_h$ carries the full four-tensor norm.
In dimension four, the Chern--Gauss--Bonnet formula gives
\begin{equation}\label{eq:flow-energy}
   \mathcal F(h)=16\pi^2\chi(N)+2e(h),
   \qquad
   e(h)=\int_N|E_h|^2\,d\mu_h .
\end{equation}
In particular, on a fixed manifold the variation of $\mathcal F$
is determined entirely by the variation of the $L^2$ norm of the
trace-free Ricci tensor.

Given an initial metric $g$ on $N$, let $h(t)$ evolve by the
negative $L^2$-gradient flow of $\mathcal F$ \cite{Streets},
\[
   \partial_t h=-\operatorname{grad}_{L^2}\mathcal F(h),
   \qquad h(0)=g.
\]
Here the gradient is taken with respect to the inner product
\[
   (A,B)_{L^2(h)}=\int_N\langle A,B\rangle_h\,d\mu_h
\]
on symmetric two-tensors. With $\Delta=\nabla^k\nabla_k$ and
the curvature convention $h^{ij}R_{ikj\ell}=\operatorname{Ric}_{k\ell}$,
the four-dimensional flow equation takes the form
\[
   \partial_t h_{ij}
      =2\Delta\operatorname{Ric}_{ij}
       -\nabla_i\nabla_jR
       +4R_{ikj\ell}\operatorname{Ric}^{k\ell}
       -R\operatorname{Ric}_{ij}
       -|E|^2h_{ij}.
\]
All quantities on the right are computed with respect to $h(t)$.

Taking the trace of the flow equation and using
$|\operatorname{Ric}|^2=|E|^2+\tfrac14R^2$, we obtain
\[
   \operatorname{tr}_h(\partial_t h)
     =2\Delta R-\Delta R
        +4|\operatorname{Ric}|^2-R^2-4|E|^2
     =\Delta R.
\]
Since $N$ is closed, the variation of the volume form gives
\[
   \frac{d}{dt}\operatorname{Vol}(N,h(t))
     =\frac12\int_N\operatorname{tr}_h(\partial_t h)\,d\mu_h
     =\frac12\int_N\Delta R\,d\mu_h=0.
\]
Thus the flow \emph{preserves the volume}. 

Along the same solution, the defining property of the gradient gives
\[
   \frac{d}{dt}\mathcal F(h(t))
     =\bigl(\operatorname{grad}_{L^2}\mathcal F(h(t)),
                         \partial_t h(t)\bigr)_{L^2(h(t))}
     =-\|\partial_t h(t)\|_{L^2(h(t))}^2.
\]
Integrating from $0$ to $t$ and using \eqref{eq:flow-energy}
at both endpoints yields
\begin{equation}\label{eq:half-dissipation}
   2e(h(t))+
   \int_0^t\|\partial_r h(r)\|_{L^2(h(r))}^2\,dr
      =2e(g).
\end{equation}
This identity provides two estimates used below. It bounds
$e(h(t))$ by its initial value, and it controls the time integral
of the squared $L^2$ speed of the flow. The latter estimate will
allow us to compare the mean scalar curvature at positive time
with that of the initial Yamabe metric.

Now we are going to run the flow on $\Sigma_j$ with initial metric being a Yamabe representative $g_j$. The following lemma ensures there is a uniform bound on a fixed time interval for all flow solutions. 
It is the central analytical tool for all proofs.

\begin{lemma}
\label{lem:positive-time-curvature}
Let $(N_j,g_j)$ be closed, connected, oriented four-manifolds with
unit-volume Yamabe metrics. Suppose that
\[
   \sup_j\chi(N_j)<\infty,\qquad
   e(g_j)\longrightarrow0,\qquad
   Y(N_j,[g_j])\longrightarrow s>0.
\]
Assume either that $s>Y_4/\sqrt2$, or that $s=Y_4/\sqrt2$ and
each $N_j$ is an integral homology four-sphere. Let $h_j(t)$ be
the $L^2$-curvature flow with $h_j(0)=g_j$. Then for sufficiently large $j$, the flow exists on $[0,1]$ and
satisfies
\begin{equation}\label{eq:half-smoothing}
   \sup_{0<t\leq1}\sqrt{t}\,
      \|\operatorname{Rm}_{h_j(t)}\|_{L^\infty(N_j,h_j(t))}
      \leq C
\end{equation}
for a constant $C$ independent of $j$.
\end{lemma}

\begin{proof}
Let $[0,T_j)$ be the maximal smooth existence interval and set
$I_j=[0,T_j)\cap[0,1]$. We first work on the intervals $I_j$.
By volume preservation, every time slice has unit volume. Moreover,
\eqref{eq:flow-energy} and \eqref{eq:half-dissipation} give
\[
   e(h_j(t))\leq e(g_j),\qquad
   \sup_j\sup_{t\in I_j}\mathcal F(h_j(t))<\infty .
\]

Let
$
   m_j(t)=\int_{N_j}R_{h_j(t)}\,d\mu_{h_j(t)}
$
denote the mean scalar curvature.
The first variation of total scalar curvature gives
\[
   m_j'(t)=
   \int_{N_j}
      \left\langle \tfrac12 R_{h_j(t)}h_j(t)
                       -\operatorname{Ric}_{h_j(t)},
                    \partial_t h_j(t)\right\rangle d\mu_{h_j(t)} .
\]
The uniform bound for $\mathcal F$ bounds the $L^2$ norm of the
first factor. Integrating in time, applying Cauchy--Schwarz, and
using \eqref{eq:half-dissipation}, we obtain
\[
   |m_j(t)-m_j(0)|
      \leq C\sqrt{t}
         \left(\int_0^t
            \|\partial_r h_j(r)\|_{L^2(h_j(r))}^2\,dr\right)^{1/2}
      \leq C\sqrt{t\,e(g_j)} .
\]
Since $g_j$ is a unit-volume Yamabe metric,
$m_j(0)=Y(N_j,[g_j])\to s$. As $t\leq1$ on $I_j$, it follows that
\begin{equation}\label{eq:half-mean}
   \sup_{t\in I_j}|m_j(t)-s|
      \leq C\sqrt{e(g_j)}+|Y(N_j,[g_j])-s|
      \longrightarrow0  \qquad \text{as $j\to \infty$}.
\end{equation}

Fix a time slice, and abbreviate $h=h_j(t)$, $m=m_j(t)$ and
$Y_h=Y(N_j,[h])$. Evaluating $\mathcal S([h])$ on a unit-volume Yamabe representative
of $[h]$ gives $24\mathcal S([h])\leq Y_h^2$. Moreover, we have
\begin{equation}\label{eq:half-scalar}
   Y_h^2\geq24\Sc([h])
     =\int_{N_j}R_h^2\,d\mu_h-12e(h)
     \geq m^2-12e(g_j).
\end{equation}
For large $j$ the right-hand side is uniformly positive. The Yamabe
constant is continuous along the smooth flow and is positive at
$t=0$, so $Y_h>0$ throughout $I_j$. In dimension four,
writing a conformal metric as $u^2h$ gives the Yamabe quotient
\[
   Q_h(u)=
   \frac{6\int_{N_j}|\nabla u|_h^2\,d\mu_h
                   +\int_{N_j}R_hu^2\,d\mu_h}
        {\bigl(\int_{N_j}|u|^4\,d\mu_h\bigr)^{1/2}}.
\]
Testing the Yamabe quotient on
the constant function gives $Y_h\leq m$. It follows that $Y_h\to s$
uniformly on $I_j$ and
\[
   \|R_h-m\|_{L^2(h)}^2
       =\int_{N_j}R_h^2\,d\mu_h-m^2
       \leq Y_h^2+12e(h)-m^2
       \leq12e(h).
\]

H\"older's inequality gives
\[
   \int_{N_j}R_hu^2\,d\mu_h
      \leq m\|u\|_{L^2(h)}^2
          +\|R_h-m\|_{L^2(h)}\|u\|_{L^4(h)}^2 .
\]
Consequently there are numbers $\delta_j\downarrow0$, after
enlarging them if necessary, such that
\begin{equation}\label{eq:half-sobolev}
   (s-\delta_j)\|u\|_{L^4(h)}^2
       \leq6\|\nabla u\|_{L^2(h)}^2
                  +(s+\delta_j)\|u\|_{L^2(h)}^2 .
\end{equation}
The same $\delta_j$ works for every $t\in I_j$ and every $u$.

Put $K_j(t)=\|\operatorname{Rm}_{h_j(t)}\|_\infty$. If
\eqref{eq:half-smoothing} fails on the existing intervals after
discarding finitely many indices, we can pass to a subsequence and
choose $t_j\in I_j$ at a running maximum of $\sqrt t\,K_j(t)$
such that
\[
   Q_j:=K_j(t_j),\qquad
   \sqrt{t_j}\,Q_j
      =\max_{0\leq t\leq t_j}\sqrt t\,K_j(t)
      \longrightarrow\infty .
\]
Choose $p_j$ where the maximum curvature at time $t_j$ is attained.
The rescaled flows
\[
   \widehat h_j(r)=Q_j h_j(t_j+r/Q_j^2),
   \qquad -1\leq r\leq0,
\]
are therefore defined for large $j$. Every original time in this
interval belongs to $[t_j/2,t_j]$, and the running-maximum property
gives
\[
   \sup_{-1\leq r\leq0}\|\operatorname{Rm}_{\widehat h_j(r)}\|_\infty
       \leq\sqrt2,\qquad
   |\operatorname{Rm}_{\widehat h_j(0)}|(p_j)=1 .
\]
Streets' smoothing estimate \cite[Theorem~1.3]{Streets} gives
uniform bounds for every spatial covariant derivative of curvature
at $r=0$.

We now fix $r=0$ and take a subsequence in $j$. Write
$\widehat h_j=\widehat h_j(0)$. Rescaling \eqref{eq:half-sobolev}
gives
\[
   (s-\delta_j)\|u\|_{L^4(\widehat h_j)}^2
     \leq6\|\nabla u\|_{L^2(\widehat h_j)}^2
       +\frac{s+\delta_j}{Q_j}\|u\|_{L^2(\widehat h_j)}^2.
\]

Fix $x\in N_j$ and write $V_j(\rho)=\operatorname{Vol}_{\widehat h_j}
(B_{\widehat h_j}(x,\rho))$. Choose a Lipschitz cutoff $\eta$ which
equals one on $B_{\widehat h_j}(x,1/2)$, vanishes outside
$B_{\widehat h_j}(x,1)$, and satisfies $0\leq\eta\leq1$ and
$|\nabla\eta|\leq2$. Applying the rescaled Sobolev inequality to
$\eta$ gives
\[
 (s-\delta_j)V_j(1/2)^{1/2}
 \leq 6\int_{N_j}|\nabla\eta|^2\,d\mu_{\widehat h_j}
      +\frac{s+\delta_j}{Q_j}\int_{N_j}\eta^2\,d\mu_{\widehat h_j}
 \leq \left(24+\frac{s+\delta_j}{Q_j}\right)V_j(1).
\]
Here the first inequality uses
$\|\eta\|_{L^4}^2\geq V_j(1/2)^{1/2}$.
Since $\delta_j\to0$, $s>0$, and $Q_j\to\infty$, it follows that
$V_j(1/2)^{1/2}\leq C V_j(1)$ uniformly in $x$ and in all
sufficiently large $j$.

The uniform curvature bound gives
$\operatorname{Ric}_{\widehat h_j}\geq-3\kappa\,\widehat h_j$
for a fixed $\kappa\geq0$. Bishop--Gromov comparison therefore yields
\[
 \frac{V_j(1)}{V_j(1/2)}
 \leq \frac{v_{-\kappa}(1)}{v_{-\kappa}(1/2)}=:C',
\]
where $v_{-\kappa}(\rho)$ is the volume of a radius-$\rho$ ball in
the simply connected four-dimensional space form of sectional
curvature $-\kappa$. Combining the two estimates gives
$V_j(1/2)^{1/2}\leq CC'V_j(1/2)$, and hence
\[
 V_j(1)\geq V_j(1/2)\geq (CC')^{-2}.
\]

The injectivity-radius estimate of Cheeger, Gromov,
and Taylor \cite{CGT} and pointed smooth compactness now give
\[
   (N_j,\widehat h_j,p_j)\longrightarrow(X,h_\infty,p)
\]
smoothly on compact sets, where $X$ is complete and
$|\operatorname{Rm}_{h_\infty}|(p)=1$.

We claim the limit is Ricci-flat. Indeed, $e$ is scale invariant in
dimension four and $e(\widehat h_j)\leq e(g_j)\to0$. In addition,
\[
   \left\|R_{\widehat h_j}-\frac{m_j(t_j)}{Q_j}\right\|_
                  {L^2(\widehat h_j)}^2
      =\|R_{h_j(t_j)}-m_j(t_j)\|_{L^2(h_j(t_j))}^2
      \leq12e(g_j).
\]
Since $m_j(t_j)/Q_j\to0$, smooth local convergence gives both
$E_{h_\infty}=0$ and $R_{h_\infty}=0$. The uniform energy bound
gives $\int_X|\operatorname{Rm}_{h_\infty}|^2\,d\mu_{h_\infty}
<\infty$. 

Passing the rescaled Sobolev inequality to the limit
on compactly supported test functions gives
\begin{equation}\label{eq:half-bubble-sobolev}
   s\|u\|_{L^4(X)}^2\leq6\|\nabla u\|_{L^2(X)}^2,
   \qquad u\in C_c^\infty(X).
\end{equation}
It follows that $X$ is noncompact, since otherwise $u\equiv1$ would
contradict \eqref{eq:half-bubble-sobolev}. A cutoff for
$B_{\rho/2}(p)\subset B_\rho(p)$ as before gives
\[
   V(\rho/2)^{1/2}\leq C\rho^{-2}V(\rho).
\]
Ricci-flatness gives $V(\rho/2)\geq 2^{-4}V(\rho)$ by
Bishop--Gromov comparison. Hence $V(\rho)\geq c\rho^4$ at every
radius. 

Therefore the limit $(X,h_\infty)$ is complete and Ricci-flat, has finite
$L^2$ curvature, and has Euclidean volume growth. By a theorem of Bando–Kasue–Nakajima
\cite[Theorem~1.5]{BKN}, it is ALE of order four. Thus, outside
a compact set, $X$ is diffeomorphic to
$(\mathbb{R}^4\setminus\overline{B}_{R_0}(0))/\Gamma$ for a finite
group $\Gamma$ acting freely on $S^3$, and the lifted metric
approaches the Euclidean metric in these coordinates. We can therefore
test the limiting Sobolev inequality on radial functions supported
in annuli far out on this end. Take radial annular cutoffs of a Euclidean Sobolev extremal.
For each fixed cutoff, dilation out along the ALE end makes its
Sobolev quotient converge to the corresponding quotient on
$\mathbb{R}^4/\Gamma$. Both integrals in the quotient are divided
by $|\Gamma|$, and the denominator is taken to the power $1/2$.
The resulting quotient is therefore $|\Gamma|^{-1/2}$ times its
Euclidean counterpart.

Hence 
\begin{equation}\label{eq:half-ale-group}
   s\leq
   \inf_{0\neq u\in C_c^\infty(X)}
       \frac{6\int_X|\nabla u|^2\,d\mu_{h_\infty}}
            {(\int_X|u|^4\,d\mu_{h_\infty})^{1/2}}
   \leq \frac{Y_4}{\sqrt{|\Gamma|}} .
\end{equation}

If $s>Y_4/\sqrt2$, then $|\Gamma|=1$. The asymptotic volume ratio
of $X$ is consequently one. Bishop--Gromov rigidity makes $X$
isometric to $\mathbb R^4$, contradicting $|\operatorname{Rm}_{h_\infty}|(p)=1$.

If $s=Y_4/\sqrt2$, the only further possibility is $|\Gamma|=2$.
A free orthogonal involution on $S^3$ is antipodal, so the link
is $P=\mathbb{RP}^3$. Fix a cross-section at a large radius.
Smooth convergence embeds this cross-section and a collar into
$N_j$ for all sufficiently large $j$. We show that this is
impossible in an integral homology four-sphere.

A two-sided embedded copy of $P$ separates an integral homology
four-sphere into two compact sides $A$ and $B$, since a
nonseparating hypersurface would have algebraic intersection one
with a closed curve. With integral coefficients, Mayer--Vietoris
gives
\[
   H_2(A)=H_2(B)=0,\qquad
   H_1(A)\oplus H_1(B)\cong H_1(P)\cong\mathbb Z/2 .
\]
After interchanging the sides, assume $H_1(A)=0$. The universal
coefficient theorem and Poincar\'e--Lefschetz duality give
$H_2(A,P)\cong H^2(A)=0$. Exactness of
\[
   H_2(A,P)\longrightarrow H_1(P)\longrightarrow H_1(A)
\]
would then force $H_1(P)=0$. A contradiction.

We have proved \eqref{eq:half-smoothing} on $I_j$ for all
sufficiently large $j$. If $T_j\leq1$ for one of these indices,
the estimate bounds curvature as $t\uparrow T_j$, since $T_j>0$.
Streets' finite-time extension criterion
\cite[Corollary~1.9]{Streets} then extends the solution past $T_j$.
Thus $T_j>1$, completing the proof.
\end{proof}

\begin{corollary}
\label{cor:fixed-time-einstein}
Under the hypotheses of Lemma~2.1, fix $t_0\in(0,1]$.
After passage to a subsequence, there exist a closed connected
four-manifold $N_\infty$, a unit-volume metric $h$ on $N_\infty$,
and diffeomorphisms $\Phi_j\colon N_\infty\to N_j$ such that
\begin{equation}\label{eq:fixed-time-limit}
    \Phi_j^*h_j(t_0)\longrightarrow h
    \qquad\text{in }C^\infty(N_\infty).
\end{equation}
Moreover, $h$ is Eintsein with scalar cuvature $R_h=s=\lim_{j\to\infty}Y(N_j,[g_j])$.
\end{corollary}

\begin{proof}
The curvature estimate of Lemma~2.1 gives a uniform curvature
bound on $[t_0/2,t_0]$. Streets' smoothing estimate
\cite[Theorem~1.3]{Streets} therefore bounds every covariant
derivative of curvature at time $t_0$, uniformly in $j$.

At this time slice, the uniform Sobolev inequality established
in the proof of Lemma~2.1 gives a lower bound for the volume of
balls of any fixed sufficiently small radius. The injectivity-radius estimate of
Cheeger, Gromov, and Taylor \cite{CGT} now applies. The diameters
are also uniformly bounded, since disjoint balls of this radius
along a minimizing geodesic each have volume at least $v_r$,
whereas the total volume is one. Smooth compactness gives
\eqref{eq:fixed-time-limit}.

The flow energy identity gives
\[
    \int_{N_\infty}|E_h|^2\,d\mu_h
      =\lim_{j\to\infty}e(h_j(t_0))
      \leq\lim_{j\to\infty}e(g_j)=0.
\]
Thus $h$ is Einstein. Its scalar curvature is constant, and the
mean scalar curvature estimate in Lemma~2.1 identifies it as
\[
    R_h=\int_{N_\infty}R_h\,d\mu_h
       =\lim_{j\to\infty}m_j(t_0)=s.
\]
\end{proof}

\begin{proof}[Proof of Theorem~A]
For the sequence $(\Sigma_j,g_j)$ constructed above, the initial
energy and scalar curvature estimates give $e(g_j)\to0$ and
$Y(\Sigma_j,[g_j])\to s_*$. Each $\Sigma_j$ is an integral homology
four-sphere, so Lemma~2.1 applies.
Corollary~\ref{cor:fixed-time-einstein}, with $t_0=1/2$, gives a
smooth limit $(N_\infty,h)$ as in \eqref{eq:fixed-time-limit}.

In the present case $N_\infty$ is an integral homology four-sphere
and $s=s_*=8\sqrt3\,\pi$. Moreover
\[
   \Sc([h])=\frac{s_*^2}{24}=8\pi^2,\qquad
   \W([h])=8\pi^2\chi(N_\infty)-\Sc([h])=8\pi^2 .
\]
Thus $[h]$ is admissible and $\beta(N_\infty,[h])=1/2$.
By \cite[Theorem~1]{CGZ}, the underlying manifold $N_\infty$
is diffeomorphic either to the standard $S^4$ or to $\mathbb{CP}^2$,
up to reversal of orientation. Since $b_2(N_\infty)=0$, only the
sphere can occur. The diffeomorphisms furnished by Corollary~2.2
therefore identify each $\Sigma_j$ in the convergent subsequence
with the standard $S^4$. Hence $\beta(\Sigma_j)=0$, contradicting
the lower bound $\beta(\Sigma_j)\geq 1/2$ in (8).
\end{proof}

\begin{proof}[Proof of Theorem~\ref{thm:main-product-gap}]
Similarly, we prove by contradiction. Suppose $(i) \Longrightarrow (iii)$ does not hold, then 
there would be admissible four-manifolds
$M_j$ with indefinite intersection form, none diffeomorphic to
$S^2\times S^2$, and classes $c_j\in\Yadm(M_j)$ such that
\[
   \frac23\leq\beta(M_j)\leq b_j:=\beta(M_j,c_j)
      <\frac23+\frac{2}{j}.
\]
Here the underlying manifolds may vary. We first determine which
characteristic numbers can occur along such a sequence. Write
$\chi_j=\chi(M_j)$, $\tau_j=\tau(M_j)$ and $Y_j=Y(M_j,c_j)$.
Evaluation on a unit-volume Yamabe metric and the round-sphere
upper bound give
\[
   8\pi^2\chi_j(1-b_j)=\Sc(c_j)
        \leq\frac{Y_j^2}{24}\leq\frac{Y_4^2}{24}=16\pi^2 .
\]
Thus $\chi_j\leq2/(1-b_j)$, whose right-hand side tends to six.
Since $\chi_j$ is integral, $\chi_j\leq6$ for large $j$.
By \eqref{eq:indef},
\[
   |\tau_j|\leq2\chi_j(b_j-2/3)\longrightarrow0 .
\]
Hence $\tau_j=0$ for large $j$. Admissibility gives $b_1(M_j)=0$,
and indefiniteness gives $b_2^+(M_j)=b_2^-(M_j)\geq1$. Passing to
a subsequence, we therefore have
\begin{equation}\label{eq:product-topology}
   \tau_j=0,\qquad \chi_j=\chi\in\{4,6\}.
\end{equation}

Let $g_j\in c_j$ be a unit-volume Yamabe metric. Since
$b_2^+(M_j)>0$ and $Y_j>0$, the modified Yamabe estimate of
\cite[Theorem~3.1(iii), equation~(3.8)]{CGZ} gives
$Y_j^2/24\leq\W^+(c_j)$. Indeed, that result supplies a
representative $\widetilde g_j\in c_j$ with
$\int R_{\widetilde g_j}^2\,d\mu_{\widetilde g_j}
\leq24\W^+(c_j)$, and testing its Yamabe quotient on constants
gives the stated inequality. Since $\tau_j=0$, the two Weyl
energies are equal. We obtain
\[
   8\pi^2\chi(1-b_j)
    =\frac{Y_j^2}{24}-\frac12e(g_j)
    \leq\frac{Y_j^2}{24}
    \leq\W^+(c_j)=4\pi^2\chi b_j .
\]
Subtracting the two endpoints and then letting $j\to\infty$ gives
\[
   0\leq e(g_j)\leq24\pi^2\chi(b_j-2/3)\longrightarrow0,
   \qquad Y_j\longrightarrow s_\chi:=8\pi\sqrt{\chi}.
\]
For both values of $\chi$ in \eqref{eq:product-topology},
\[
   \frac{s_\chi}{Y_4}=\sqrt{\frac{\chi}{6}}>\frac1{\sqrt2}.
\]
Thus the strict-threshold case of
Lemma~\ref{lem:positive-time-curvature} applies.

Applying Corollary~\ref{cor:fixed-time-einstein} at $t_0=1/2$,
we obtain a closed unit-volume Einstein manifold $(M_\infty,h)$
with $R_h=s_\chi$. After passage to a subsequence, $M_\infty$ is
diffeomorphic to each $M_j$ in the subsequence.
 The intersection form of
$M_\infty$ is indefinite and $\chi(M_\infty)=\chi$. Evaluating
\eqref{eq:characteristic} on $h$ gives
\[
   \Sc([h])=\frac{s_\chi^2}{24}
      =\frac{8\pi^2\chi}{3},\qquad
   \beta(M_\infty,[h])
      =1-\frac{\Sc([h])}{8\pi^2\chi}
      =\frac23 .
\]
By \cite[Theorem~C]{CGZ}, $(M_{\infty}, [h])$ is conformal equivalent to 
$(S^2\times S^2, [g_{rd}\oplus g_{rd}])$. Then for all sufficiently large $j$, $M_j$ is diffeomorphic to $S^2\times S^2$, a contradiction.
\end{proof}

\section{Strict lower bounds for the one- and two-point blow-ups
of $\CP^2$}
\label{sec:delpezzo-gaps}

Let $X_k=\CP^2\#k\overline{\CP}^{\,2}$, equipped with its complex
orientation. We have $b_2^+(X_k)=1$,
$\chi(X_k)=3+k$, and $\tau(X_k)=1-k$. Gursky's inequality
\eqref{eq:lower} therefore gives
\[
   \beta(X_1)\geq\frac23,\qquad
   \beta(X_2)\geq\frac{23}{30}.
\]
We show that both inequalities are strict, even after taking the
infimum over all admissible conformal classes. For $X_1$, this
follows from Theorem~B and the parity of the intersection form.
For $X_2$, we show that a sequence approaching $23/30$ would
produce a positive K\"ahler--Einstein metric on a del Pezzo
surface of degree seven. The Page and Chen--LeBrun--Weber classes
give explicit upper bounds for the two infima.

\begin{theoremx}
\label{thm:main-blowups}
Let $g_{\mathrm{P}}$ and $g_{\mathrm{CLW}}$ denote the Page metric
on $X_1$ and the Chen--LeBrun--Weber metric on $X_2$, respectively.
Then
\begin{equation}\label{eq:intervals}
   \begin{aligned}
      \frac23
         &<\beta(X_1)\leq\beta(X_1,[g_{\mathrm{P}}])\approx0.682,\\
      \frac{23}{30}
         &<\beta(X_2)\leq\beta(X_2,[g_{\mathrm{CLW}}])\approx0.776.
   \end{aligned}
\end{equation}
\end{theoremx}

\begin{proof}

\emph{On $X_1$.}

The integral intersection form of $X_1$ is
$\langle 1\rangle\oplus\langle-1\rangle$, which is odd. The
intersection form of $S^2\times S^2$ is even. Thus $X_1$ is
indefinite and is not diffeomorphic to $S^2\times S^2$.
Theorem~\ref{thm:main-product-gap} excludes the entire interval
$[2/3,2/3+\varepsilon_2)$ and gives
\[
   \beta(X_1)\geq\frac23+\varepsilon_2.
\]

\medskip

\emph{On $X_2$.}

Suppose that $\beta(X_2)=23/30$. Choose admissible classes $c_j$ 
on $X_2$ with $b_j:=\beta(X_2,c_j)\to23/30$, and let $g_j\in c_j$
be unit-volume Yamabe metrics. Write
\[
   Y_j=Y(X_2,c_j)=R_{g_j},\qquad
   e(g_j)=\int_{X_2}|E_{g_j}|^2\,d\mu_{g_j}.
\]
We first determine the limiting Yamabe constant and the $L^2$ norm
of the trace-free Ricci curvature. Since $\chi(X_2)=5$ and
$\tau(X_2)=-1$, the characteristic number formulas
\eqref{eq:characteristic} give
\[
   \W^+(c_j)=20\pi^2 b_j-6\pi^2,\qquad
   \Sc(c_j)=40\pi^2(1-b_j).
\]
Set
\[
   L=\frac{4\pi^2}{3}\bigl(2\chi(X_2)+3\tau(X_2)\bigr)
     =\frac{28\pi^2}{3},\qquad
   a_j=20\pi^2\left(b_j-\frac{23}{30}\right).
\]
Then $a_j\geq0$, $a_j\to0$, and the preceding identities become
\begin{equation}\label{eq:x2-defect}
   \W^+(c_j)=L+a_j,\qquad \Sc(c_j)=L-2a_j.
\end{equation}

Because $b_2^+(X_2)>0$ and $Y_j>0$, the modified Yamabe estimate
used in the proof of Theorem~\ref{thm:main-product-gap} gives
$Y_j^2/24\leq\W^+(c_j)$
\cite[Theorem~3.1(iii), equation~(3.8)]{CGZ}.
Evaluating $\Sc(c_j)$ on the ordinary Yamabe representative $g_j$
and using \eqref{eq:x2-defect}, we obtain
\[
   L-2a_j=\frac{Y_j^2}{24}-\frac12e(g_j)
       \leq\frac{Y_j^2}{24}\leq L+a_j.
\]
It follows that
\begin{equation}\label{eq:x2-yamabe-limit}
   e(g_j)\leq6a_j\longrightarrow0,\qquad
   Y_j\longrightarrow s:=\sqrt{24L}=4\sqrt{14}\,\pi .
\end{equation}
The comparison with the bubbling threshold is strict:
\[
   \frac{s}{Y_4}=\sqrt{\frac{7}{12}}>\frac1{\sqrt2}.
\]

The strict inequality $s>Y_4/\sqrt{2}$ places this sequence in
the first case of Lemma~2.1. By
Corollary~\ref{cor:fixed-time-einstein}, the time-$1/2$ slices
have a subsequence converging smoothly, up to diffeomorphism,
to a unit-volume Einstein metric $h$ with $R_h=s$. We use one
of the convergence diffeomorphisms to transport $h$ to $X_2$,
which we henceforth equip with its complex orientation.

To identify the conformal class of the limit, we recover its
self-dual Weyl energy from the characteristic numbers. Since $h$
is Einstein and has unit volume,
\[
   \Sc([h])=\frac{R_h^2}{24}=L.
\]
Adding the two identities in \eqref{eq:characteristic} gives
\[
   \Sc([h])+2\W^+([h])
      =4\pi^2\bigl(2\chi(X_2)+3\tau(X_2)\bigr)=3L.
\]
Hence $\W^+([h])=L$. Thus $[h]$ attains equality in Gursky's inequality. By its
equality characterization \cite{Gursky}, there is a complex
structure $J$ inducing the chosen orientation on $X_2$ for which
$[h]$ contains a positive K\"ahler--Einstein metric. To identify
this complex surface, observe that positivity of the Ricci form
gives $c_1(X_2,J)>0$. Hence $(X_2,J)$ is a del Pezzo surface,
and its degree is
\[
    c_1(X_2,J)^2
       =2\chi(X_2)+3\tau(X_2)
       =2\cdot5+3(-1)=7.
\]
By the classification of smooth del Pezzo surfaces
\cite[\S1]{LeBrunWeyl}, it is the blow-up of $\mathbb{CP}^2$
at two distinct points. All such pairs of points are equivalent
under projective automorphisms, so the resulting complex surface
is unique up to biholomorphism. Its holomorphic automorphism
group has a nonreductive identity component. Matsushima's theorem
therefore excludes a positive K\"ahler--Einstein metric on
$(X_2,J)$ \cite[\S1]{CLW}. This contradiction proves
$\beta(X_2)>23/30$.

The Page and Chen--LeBrun--Weber metrics are positive Einstein
metrics \cite{Page,CLW}. Each has positive Yamabe constant and
$\Sc=\tfrac1{24}\int R^2\,d\mu>0$, so its conformal class is
admissible. The latter inequalities in \eqref{eq:intervals} follow
from the definition of $\beta(X_k)$.

For the numerical values, let $g_K$ be the extremal K\"ahler
representative of either Einstein conformal class and set
\[
   A(g_K)=\frac{1}{32\pi^2}\int_{X_k}R_{g_K}^2\,d\mu_{g_K}.
\]
The K\"ahler surface identity $\|W^+_{g_K}\|^2=R_{g_K}^2/24$
gives $\W^+([g_K])=4\pi^2 A(g_K)/3$. The signature formula
therefore yields
\[
   \beta(X_k,[g_K])
     =\frac{A(g_K)}{3\chi(X_k)}
        -\frac{3\tau(X_k)}{2\chi(X_k)} .
\]
For the Page class, set $k=1$ in the Hirzebruch-surface formula
of \cite[\S7, equation~(7.1)]{LeBrunCalabi}. In the parameter
$\alpha$ used there, the virtual action is
\[
    A_1(\alpha)
    =\frac{2\alpha^3+7\alpha^2+8\alpha+3/2}
           {\alpha^2+\alpha+1/6}.
\]
Its positive critical point is
$\alpha_*\approx0.458$, giving
$A_1(\alpha_*)\approx8.179$. For the
Chen--LeBrun--Weber class, the first positive critical point
of the function $f$ in \cite[\S3, Lemma~3]{CLW} is
$x_*\approx0.9577128052$, with $f(x_*)\approx7.136474471$.
Substituting into the preceding formula gives
\[
\begin{aligned}
    \beta(X_1,[g_{\mathrm{P}}])
       &=\frac{A_1(\alpha_*)}{12}
         \approx0.682,\\
    \beta(X_2,[g_{\mathrm{CLW}}])
       &=\frac{3}{10}+\frac{f(x_*)}{15}
         \approx0.776.
\end{aligned}
\]

\end{proof}


\begin{thebibliography}{99}
\small
















\bibitem{BKN}
S. Bando, A. Kasue, and H. Nakajima,
\emph{On a construction of coordinates at infinity on manifolds
with fast curvature decay and maximal volume growth},
Invent. Math. \textbf{97} (1989), 313--349.

\bibitem{BroderLafuente}
K. Broder and R. A. Lafuente,
\emph{$T^2$-invariant Einstein metrics on closed $4$-manifolds},
preprint, arXiv: 2610.18219.


\bibitem{CGYsphere}
S.-Y. A. Chang, M. J. Gursky, and P. C. Yang,
\emph{A conformally invariant sphere theorem in four dimensions},
Publ. Math. Inst. Hautes \'Etudes Sci. \textbf{98} (2003), 105--143.


\bibitem{CGZ}
S.-Y. A. Chang, M. Gursky, and S. Zhang,
\emph{A conformally invariant gap theorem characterizing $\CP^2$ via the Ricci flow},
Math. Z. \textbf{294} (2020), 721--746.


\bibitem{CQY}
S.-Y. A. Chang, J. Qing, and P. Yang,
\emph{On a conformal gap and finiteness theorem for a class of four-manifolds},
Geom. Funct. Anal. \textbf{17} (2007), 404--434.

\bibitem{CGT}
J. Cheeger, M. Gromov, and M. Taylor,
\emph{Finite propagation speed, kernel estimates for functions of
the Laplace operator, and the geometry of complete Riemannian
manifolds},
J. Differential Geom. \textbf{17} (1982), 15--53.


\bibitem{ChenZhu}
B.-L. Chen and X.-P. Zhu,
\emph{A conformally invariant classification theorem in four
dimensions},
Comm. Anal. Geom. \textbf{22} (2014), 811--831.


\bibitem{CLW}
X. Chen, C. LeBrun, and B. Weber,
\emph{On conformally K\"ahler, Einstein manifolds},
J. Amer. Math. Soc. \textbf{21} (2008), 1137--1168.

\bibitem{ChuStern}
A. C.-P. Chu and D. Stern,
\emph{A variational construction for Einstein metrics on connected
sums of $S^2\times S^2$},
preprint, arXiv: 2610.29882.

\bibitem{Donaldson}
S. K. Donaldson,
\emph{An application of gauge theory to four-dimensional topology},
J. Differential Geom. \textbf{18} (1983), 279--315.

\bibitem{Freedman}
M. H. Freedman,
\emph{The topology of four-dimensional manifolds},
J. Differential Geom. \textbf{17} (1982), 357--453.


\bibitem{Gursky}
M. J. Gursky,
\emph{The Weyl functional, de Rham cohomology and K\"ahler--Einstein metrics},
Ann. of Math. (2) \textbf{148} (1998), 315--337.



\bibitem{Koca}
C. Koca,
\emph{On Conformal Geometry of K\"ahler Surfaces},
Ph.D. thesis, Stony Brook University, 2012.


\bibitem{LeBrunCalabi}
C. LeBrun,
\emph{Calabi energies of extremal toric surfaces},
Surv. Differ. Geom. \textbf{18} (2013), 195--226.


\bibitem{LeBrunWeyl}
C. LeBrun,
\emph{Weyl curvature, del Pezzo surfaces, and almost-K\"ahler geometry},
J. Geom. Anal. \textbf{25} (2015), 1744--1772.


\bibitem{LNN}
C. LeBrun, S. Nayatani, and T. Nitta,
\emph{Self-dual manifolds with positive Ricci curvature},
Math. Z. \textbf{224} (1997), 49--63.





\bibitem{LiSun}
M. Li and S. Sun,
\emph{Compact toric Einstein four-manifolds},
preprint, arXiv: 2610.10425.

\bibitem{MorganTian}
J. Morgan and G. Tian,
\emph{The Geometrization Conjecture},
Clay Mathematics Monographs, vol.~5,
American Mathematical Society and Clay Mathematics Institute, 2014.


\bibitem{Ozuch}
T. Ozuch,
\emph{Positive Einstein metrics on infinitely many four-manifolds},
preprint, arXiv: 2610.08052.

\bibitem{Page}
D. N. Page,
\emph{A compact rotating gravitational instanton},
Phys. Lett. B \textbf{79} (1978), 235--238.

\bibitem{Streets}
J. Streets,
\emph{The long time behavior of fourth order curvature flows},
Calc. Var. Partial Differential Equations \textbf{46} (2013),
39--54.


\bibitem{Tian}
G. Tian,
\emph{On Calabi's conjecture for complex surfaces with positive first Chern class},
Invent. Math. \textbf{101} (1990), 101--172.


\bibitem{TianYau}
G. Tian and S.-T. Yau,
\emph{K\"ahler--Einstein metrics on complex surfaces with $C_1>0$},
Comm. Math. Phys. \textbf{112} (1987), 175--203.



\bibitem{Zhang}
S. Zhang,
\emph{Some conformally invariant gap theorems for Bach-flat 4-manifolds},
Calc. Var. Partial Differential Equations \textbf{60} (2021), Paper No.~151.




















\end{thebibliography}
\end{document}